\documentclass[aps,prl,twocolumn,showpacs,superscriptaddress,groupedaddress]{revtex4}  % for review and submission
\usepackage{graphicx}  % needed for figures
\usepackage{dcolumn}   % needed for some tables
\usepackage{bm}        % for math
\usepackage{amssymb}
\usepackage{overpic}
\usepackage{graphics}
\usepackage{epsfig}
\usepackage{amsmath}
\usepackage{amsfonts}
\usepackage{paralist}

\newtheorem{fact}{Fact}

\begin{document}

\title{Isospectral Reductions of Dynamical Networks}

%    author one information
\author{L. A. Bunimovich}
\address{ABC Math Program and School of Mathematics, Georgia Institute of Technology, Atlanta, USA}
\email{bunimovh@math.gatech.edu}

%    author one information
\author{B. Z. Webb}
\address{ABC Math Program and School of Mathematics, Georgia Institute of Technology, Atlanta, USA}
\email{bwebb@math.gatech.edu}

\begin{abstract}
We present a general and flexible procedure which allows for the reduction (or expansion) of any dynamical network while preserving the spectrum of the network's adjacency matrix. Computationally, this process is simple and easily implemented for the analysis of any network. Moreover, it is possible to isospectrally reduce a network with respect to any network characteristic including centrality, betweenness, etc. This procedure also establishes new equivalence relations which partition all dynamical networks into spectrally equivalent classes. Here, we present general facts regarding isospectral network transformations which we then demonstrate in simple examples. Overall, our procedure introduces new possibilities for the analysis of networks in ways that are easily visualized.
\end{abstract}


\maketitle

Many, if not the majority of networks that occur in nature are very large in the sense that they contain many nodes and many edges \cite{1,5,4}. A considerable number of these networks (neural networks, metabolic networks, power stations, etc.) are also dynamic, i.e. each node has an associated state that changes in time according to its intrinsic evolution and interactions with other nodes. 

Because of the complex structure of most real-world networks there is a strong motivation for developing methods that simplify the topology of a network while maintaining some of the network's basic characteristics. In this regard, perhaps, the most fundamental characteristic of a dynamical network is the spectrum of its adjacency matrix \cite{2,8}. In a dynamical network each edge of the graph of interactions (topology of a network) can be assigned a weight related to the strength of the corresponding interaction. In practice, such edge weights typically represent the linearization of the associated interaction. These weights make up network's weighted adjacency matrix which forms the backbone of any dynamical network. An important question then is how and to what extent can a network be simplified (reduced in number of nodes and edges) while preserving the spectrum (eigenvalues) of its adjacency matrix.

At first it may seem that both simplifying (i.e. reducing) a network and maintaining its spectrum is not possible. One obvious issue is dimensionality. Indeed, a network with $n$ nodes has an $n\times n$ adjacency matrix so reducing a network to $m<n$ nodes would seem to necessarily reduce the number of eigenvalues associated with the network. However, despite this and other possible complications we will demonstrate that it is in fact possible to do just this. 

The process of \textit{isospectral reduction}, which we propose, moreover allows for the introduction of new equivalence relations between dynamical networks where two networks are equivalent if they can be isospectrally reduced to the same network. This procedure is also quite flexible as it is possible to reduce a network to a smaller network on any nonempty subset of its original set of nodes. Our procedure, therefore, allows one to reduce a network with regard to characteristics such as centrality, betweenness, clustering, etc. \cite{4} e.g. over nodes with maximal centrality. In this respect, the method of isospectral reductions allows experts (e.g. physicists, biologists, chemists, engineers, etc.) to reduce (simplify) the dynamical networks with which they are concerned relative to any network characteristic they find interesting while preserving the network's major dynamical characteristics. 

Additionally, this approach to network analysis is visually informative as one can directly see and compare various reductions of the same network or reductions of different networks. Also, if one wishes to reduce edges instead of nodes one need only interchange nodes and edges and perform the same operation.

As a dual procedure we similarly demonstrate that one can enlarge (rather than reduce) a dynamical network, thereby making it more sparse, while again maintaining its spectrum up to a set known in advance. Here, as in the case of isospectral reductions, this collection of eigenvalues is immediately known from the topology and weights of the network.

Our results on isospectral transformations are characterized by two important features. First, these results can be rigorously proven (see \cite{3}). Therefore, we refer to them here as Facts. Second, our procedure is algorithmic and importantly easy to implement. Therefore, the major goal of this paper is to present the network community with a new and flexible tool for the analysis of dynamical networks which could be immediately put to use in the study of real networks. The only information needed to carry out this procedure is the topology of a specific network.

\begin{figure*}
    \begin{overpic}[scale=.37,angle=270]{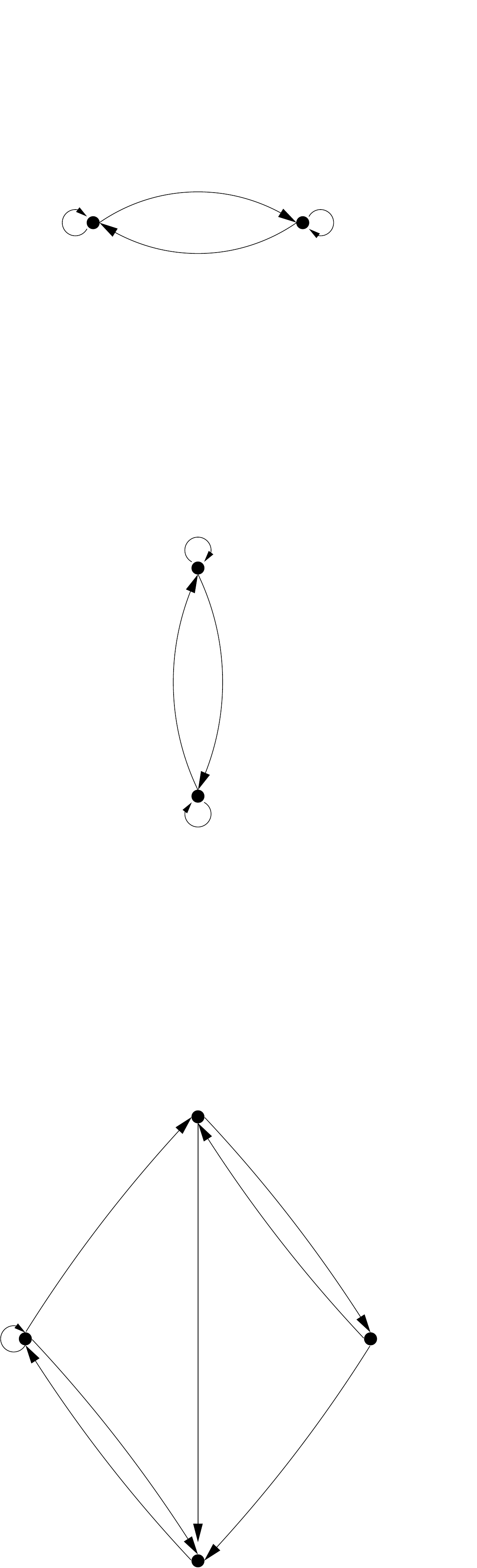}
    \put(13.5,3){$G$}
    \put(4,0){$\sigma(G)=\{2,-1,0,0\}$}
    \put(-2,18){$v_1$}
    \put(29.5,18){$v_2$}
    \put(16.5,28.5){$v_3$}
    \put(11,6.5){$v_4$}
    \put(5,25){$\frac{1}{\lambda}$}
    \put(9.5,21.5){$\lambda$}
    \put(12,29){$1$}
    \put(23.25,24.5){$2$}
    \put(14,19){$\frac{\lambda}{2}$}
    \put(3,10){$-\frac{\lambda}{6}$}
    \put(22,9.75){$-3$}
    \put(16.75,13.25){$-\frac{1}{3}$}

    \put(43,17.5){$\frac{1}{\lambda-1}$}
    \put(54.5,20.5){$\frac{2\lambda}{\lambda-1}$}
    \put(53.5,14){$\frac{1}{2}+\frac{\lambda}{2}$}
    \put(53,10){$\mathcal{R}_{S}(G)$}
    \put(66,17.5){$\frac{1}{\lambda}$}
    \put(46,7){$\sigma(\mathcal{R}_{S}(G))=\{2,-1\}$}
    \put(48.5,19){$v_1$}
    \put(62,19){$v_2$}

    \put(82,3){$\mathcal{R}_{T}(G)$}
    \put(74,0){$\sigma(\mathcal{R}_{T}(G))=\{2,-1,0\}$}
    \put(87,24){$v_3$}
    \put(87,10.5){$v_4$}
    \put(85,7){$\frac{1}{\lambda}$}
    \put(73.5,17.5){$-\frac{1}{6\lambda}-\frac{1}{6\lambda^2}$}
    \put(88.5,17.5){$-\frac{6}{\lambda}$}
    \put(81,28){$1+\frac{1}{\lambda}+\frac{1}{\lambda^2}$}
    \end{overpic} 
\caption{The graph $G$ (left) and its reductions over $S=\{v_1,v_3\}$ (center) and $T=\{v_2,v_4\}$ (right) where each spectrum is indicated. Here all edge weights of $G$ are equal to 1.}
\end{figure*}
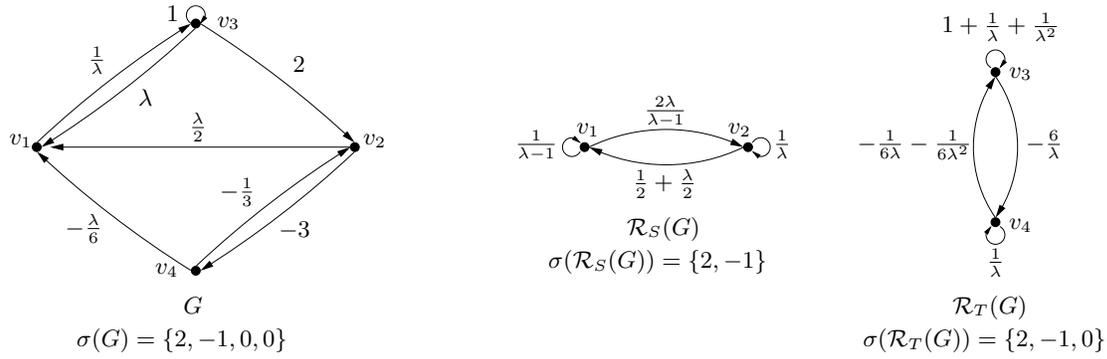 

In this paper, every network is described by its set of nodes (vertices) $V$, directed edges $E$, and weights $\omega(e)$ of the edges $e$. Thus, a dynamical network is identified with a weighted and directed graph $G=(V,E,\omega)$. An edge from vertex $v_i$ to $v_j$ is denoted by $e_{ij}$. The \textit{weighted adjacency matrix} of $G$ is given entrywise by $M(G)_{ij}=\omega(e_{ij})$. The \textit{spectum} $\sigma(G)$, or the collection of eigenvalues of $M(G)$ are the solutions to the equation
\begin{equation}\label{eq1}
\det(M(G)-\lambda I)=0.
\end{equation}
Our procedure allows edge weights of a network to be maintained under reduction if all are equal to 1 (i.e. if a network is unweighted) or weights to be positive numbers, etc. However, in the most general case edge weights are allowed to be functions of the spectral parameter $\lambda$ in equation (\ref{eq1}). Various examples will be presented below.

Of primary importance is the fact that a typical network $G$ cannot be reduced over an arbitrary subset of its vertices by means of a single reduction. Instead, a network $G$ can be reduced over any arbitrary subset of its vertices via a sequence of isospectral reductions where each reduction is achieved by reducing over a particular subset of the network's vertices. Such subsets are called structural sets and are described as follows.

First, we remove all loops (closed paths containing a single vertex) from $G$ and denote the remaining graph by $\ell(G)$. A subset $S\subseteq V$ is a \textit{structural set} of $G$ if there are no closed paths in $\ell(G)$ on $\bar{S}=V-S$ and $\omega(e_{ii})\neq\lambda$ for each $v_i\in \bar{S}$. For example, in Fig. 1 the vertex sets $S=\{v_1,v_2\}$ and $T=\{v_3,v_4\}$ are structural sets of $G$ whereas $\{v_1,v_3\}$ is not. From a computational point of view it is straightforward to determine, for a given network, which subsets of its vertices are structural sets and which are not. The formal procedure (algorithm) for isospectrally reducing a graph is as follows.

A \textit{path} in a graph $G=(V,E,\omega)$ is an ordered sequence of distinct vertices $v_1,\dots,v_m$ such that $e_{i,i+1}$ are edges in $E$ for each $1\leq i\leq m-1$. The \textit{interior vertices} of a path are the path's vertices except the first and last. For a structural set $S=\{v_1,\dots,v_n\}$ let $\mathcal{B}_{ij}(G;S)$ be the collection of all paths from $v_i$ to $v_j$ for $1\leq i,j\leq n$ which contain no interior vertices in $S$. Now consider the union of such paths over all pairs $v_i$ and $v_j$ given by
\begin{equation}\label{eq2}
\mathcal{B}_S(G)=\bigcup_{1\leq i,j\leq n}\mathcal{B}_{ij}(G;S)
\end{equation}
An \textit{isospectral reduction} of $G=(V,E,\omega)$ over the structural set $S=\{v_1,\dots,v_n\}$ is the network $\mathcal{R}_S(G)=(S,\mathcal{E},\mu)$, where there is an edge $e_{ij}\in\mathcal{E}$ between $v_i$ and $v_j$ of $S$ if and only if there is at least one path in $\mathcal{B}_{ij}(G;S)$. If $e_{ij}\in\mathcal{E}$ then $e_{ij}$ has weight equal to
\begin{equation}\label{eq3}
\mu(e_{ij})=\sum_{\beta\in\mathcal{B}_{ij}(G;S)}\omega(e_{12})\prod_{i=0}^{m-1}\frac{\omega(e_{i,i+1})}{\lambda-\omega(e_{ii})}
\end{equation}
where the sum is taken over all paths $\beta=v_1,\dots,v_m$.

\begin{figure*}
    \begin{overpic}[scale=.34,angle=270]{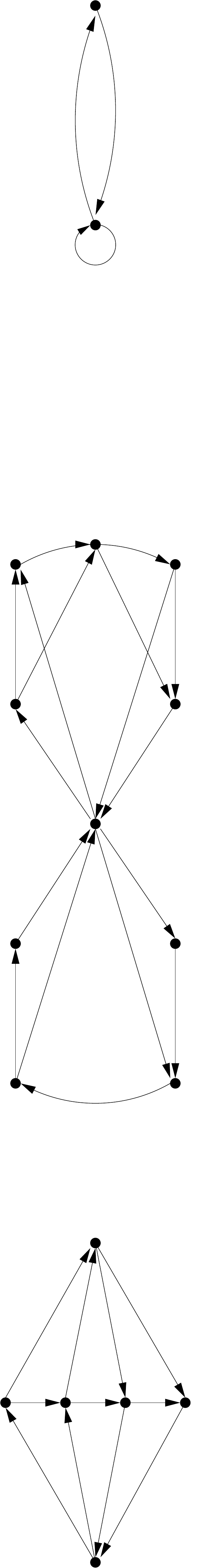}
    \put(9.5,-1.5){$H$}
    \put(44.5,-1.5){$H^\prime$}
    \put(-2,6.5){$v_1$}
    \put(21.5,6.5){$v_2$}
    \put(47,8){$v_1$}
    \put(66,6.5){$v_2$}
    \put(86,8){$v_1$}
    \put(98.5,8){$v_2$}
    \put(90,9){$\frac{2}{\lambda}+\frac{1}{\lambda^2}$}
    \put(90,3.5){$\frac{2}{\lambda}+\frac{1}{\lambda^2}$}
    \put(71.5,6.5){$\frac{1}{\lambda^2}+\frac{2}{\lambda^3}+\frac{1}{\lambda^4}$}
    \put(76,-1){$\mathcal{R}_{\{v_1,v_2\}}(H)=\mathcal{R}_{\{v_1,v_2\}}(H^\prime)$}
    \end{overpic} 
\caption{As $\mathcal{R}_{\{v_1,v_2\}}(H)=\mathcal{R}_{\{v_1,v_2\}}(H^\prime)$ the networks $H$ and $H^\prime$ are spectrally equivalent. Here all edge weights of $H$ and $H^\prime$ are equal to 1 and only the vertices of the structural set $\{v_1,v_2\}$ are labeled.}
\end{figure*}
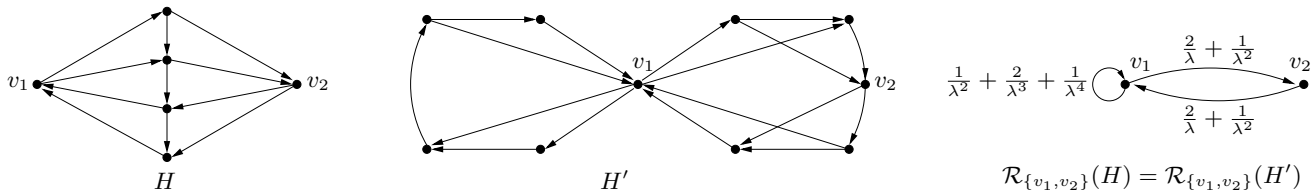

\begin{fact} \textbf{(Preservation of Spectral Information)}
The spectrum of a network $G=(V,E,\omega)$ is preserved when $G$ is reduced over the structural set $S$ up to the set $\sigma(G,S)$ which is the collection of weights of loops of the vertices in $\bar{S}$. Therefore, this set is known in advance and can be directly read from the initial network (see theorem 3.4 in \cite{3}).
\end{fact}

For example, in Fig. 1 consider the structural set $S=\{v_1,v_2\}$. Here, one can see that $\sigma(G,S)=\{0,0\}$ which is the difference between $\sigma(G)$ and $\sigma(\mathcal{R}_S(G))$.

With the concept of an isospectral reduction in place it is natural to define an \textit{isospectral expansion} of a graph $G$ as a graph $H$ where $\mathcal{R}_T(H)=G$ for some structural set $T$ of $H$. Such expansions can be carried out by expanding edges into paths or multiple paths with the correct product and sum given in (\ref{eq3}).

Importantly, as a reduced network $\mathcal{R}_S(G)$ will have its own structural sets, it possible to sequentially reduce a network. However, given that different sequences of reductions seemingly lead to different reductions, a natural question then is whether one can have some sort of control of the resulting reduced network in this huge variety of possible reductions. That is, by which sequence of reductions should one reduce a network to a smaller network on a particular vertex set that one finds interesting? As it happens, sequences of reductions have the following remarkable property which resolves this key issue.

\begin{fact} \textbf{(Commutativity of Sequential Network Reductions)} In a sequence of reductions order and number of reductions do not matter. What matters is only the final subset of vertices over which the network is reduced (see theorem 3.5 in \cite{3}). 
\end{fact}

More formally this fact could be expressed in the following way. Let the sets $S_m\subset S_{m-1}\subset\dots\subset S_1\subset V$ be structural sets in a sequence of reductions on the graph $G=(V,E,\omega)$. Then $\mathcal{R}(G;S_1,\dots,S_m)=\mathcal{R}_{S_m}(G)$. That is, regardless of the specific sequence of reductions over which the graph $G$ is reduced to a graph on the vertices $S_m$ the result is always the same. 

Besides allowing one to simplify (reduce) a network, our procedure can also be used to establish new equivalence relations between different networks. Namely, we say two networks are \textit{spectrally equivalent} if they can be isospectrally reduced to the same network. As an example, in Fig. 2 we show the two spectrally equivalent networks $H$ and $H^\prime$. However, note that even in this relatively simple example it is rather improbable that one would recognize this fact if it were not pointed out. As real-world networks are much more complicated the following fact is quite remarkable and instructive on how to deal with networks of any type. 

\begin{fact} \textbf{(Network Equivalence)}
Choose a rule $\tau$ that uniquely selects a subset of nodes $\tau(G)\subseteq V$ of any network $G=(V,E,\omega)$. Then the rule $\tau$ induces an equivalence relation on the space of all networks in the sense that two networks $G$ and $H$ belong to the same equivalence class if they are spectrally equivalent under $\tau$ i.e. $\mathcal{R}_{\tau(G)}(G)=\mathcal{R}_{\tau(H)}(H)$ (see theorem 3.9 in \cite{3}).
\end{fact}

For instance, the choice of all vertices of a network with minimal out degree is a rule that selects a unique subset of nodes of any network. Conversely, the choice of an arbitrary vertex of a network is not. Here, it is again the role of the expert to choose a specific rule for analysis and comparison of concrete networks observing that the flexibility of our procedure allows for many such rules. For a natural starting point one could use well known network characteristics (e.g. degree, centrality, betweenness) to develop such rules. The only requirement is that such rules must single out a unique subset of nodes. 

So far in this paper we have considered network transformations in which the network's edge weights become more complicated as the network is reduced. It could be argued from this point of view that our procedure simply trades network complexity for the complexity of these edge weights. However, this is not the case. 

\begin{fact} \textbf{(Isospectral Reductions Over Fixed Weight Sets)}
It is usually possible to isospectrally reduce a network while maintaining the network's set of edge weights. Moreover, only the number of zeros in the spectrum of $G$ can change under such transformations. (see theorem 4.5 in \cite{3})
\end{fact}
The procedure for reducing a network $G=(V,E,\omega)$ over a fixed weight set is as follows. First, let $S=\{v_1,\dots,v_n\}$ be a structural set of $G$ with the restriction that each vertex of $\bar{S}$ has no loop. 

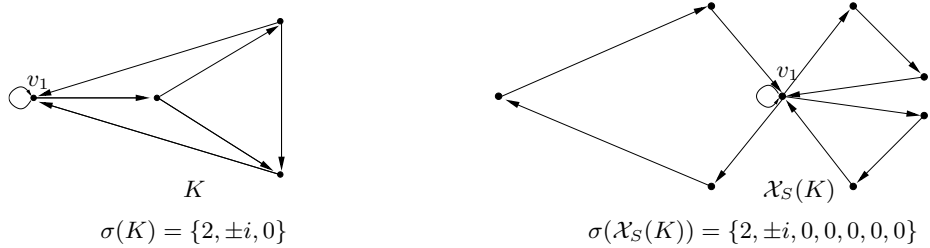
\begin{figure*}
    \begin{overpic}[scale=.27,angle=270]{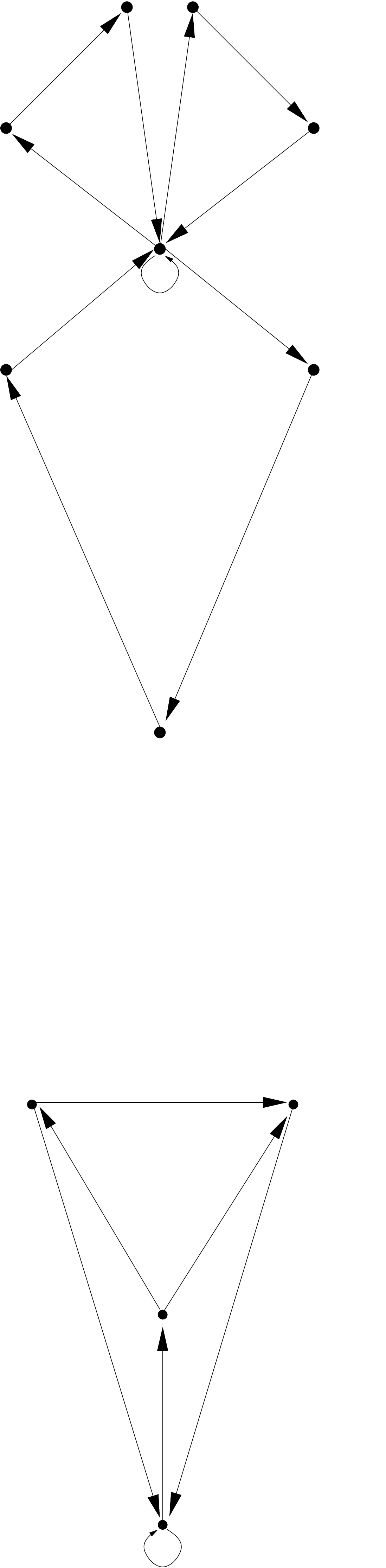}
    \put(19,3.5){$K$}
    \put(10,-1){$\sigma(K)=\{2,\pm i,0\}$}
    \put(82,3.5){$\mathcal{X}_S(K)$}
    \put(63,-1){$\sigma(\mathcal{X}_S(K))=\{2,\pm i,0,0,0,0,0\}$}
    \put(2,15.5){$v_1$}
    \put(83.5,16.5){$v_1$}
    \end{overpic} 
\caption{For $S=\{v_1\}$, $\mathcal{X}_S(K)$ is an expansion (sparsification) of $K$ over the weight set $\{1\}$, i.e. the edge weight 1 of each edge is maintained. Here only the structural set $S=\{v_1\}$ is labeled.}
\end{figure*}

\textit{Step 1:} For each path $P=v_1,\dots,v_m$ in $\mathcal{B}_S(G)$ let $\pi(P)$ be the product of the edge weights along this path. Then reweight the edges of this path by giving the first edge of $P$ weight $\pi(P)$ and the last $m-2$ edges weight $1$. Denote the set of reweighted paths terminating at the vertex $v_j$ by 
\begin{equation}\label{eq4}
\pi\mathcal{B}_j(G;S)=\bigcup_{1\leq i\leq n}\pi\mathcal{B}_{ij}(G;S)
\end{equation}
where $\pi\mathcal{B}_{ij}(G;S)$ are the reweighted paths from $v_i$ to $v_j$.

\textit{Step 2:} For $1\leq j\leq n$, combine the last interior vertex of each path in $\pi\mathcal{B}_j(G;S)$ to a single vertex leaving a single edge between this vertex and $v_j$ of weight 1. Similarly, combine the second to last interior vertex of each path in this manner and so on until all interior vertices of paths in $\pi\mathcal{B}_j(G;S)$ have been combined in this way.

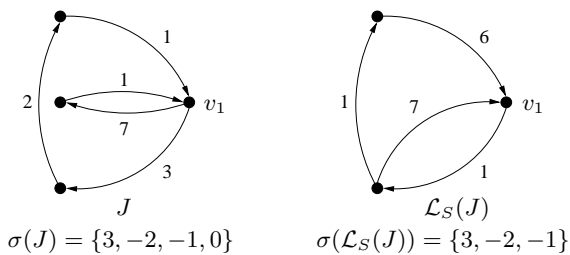
\begin{figure}[b]
    \begin{overpic}[scale=.45,angle=270]{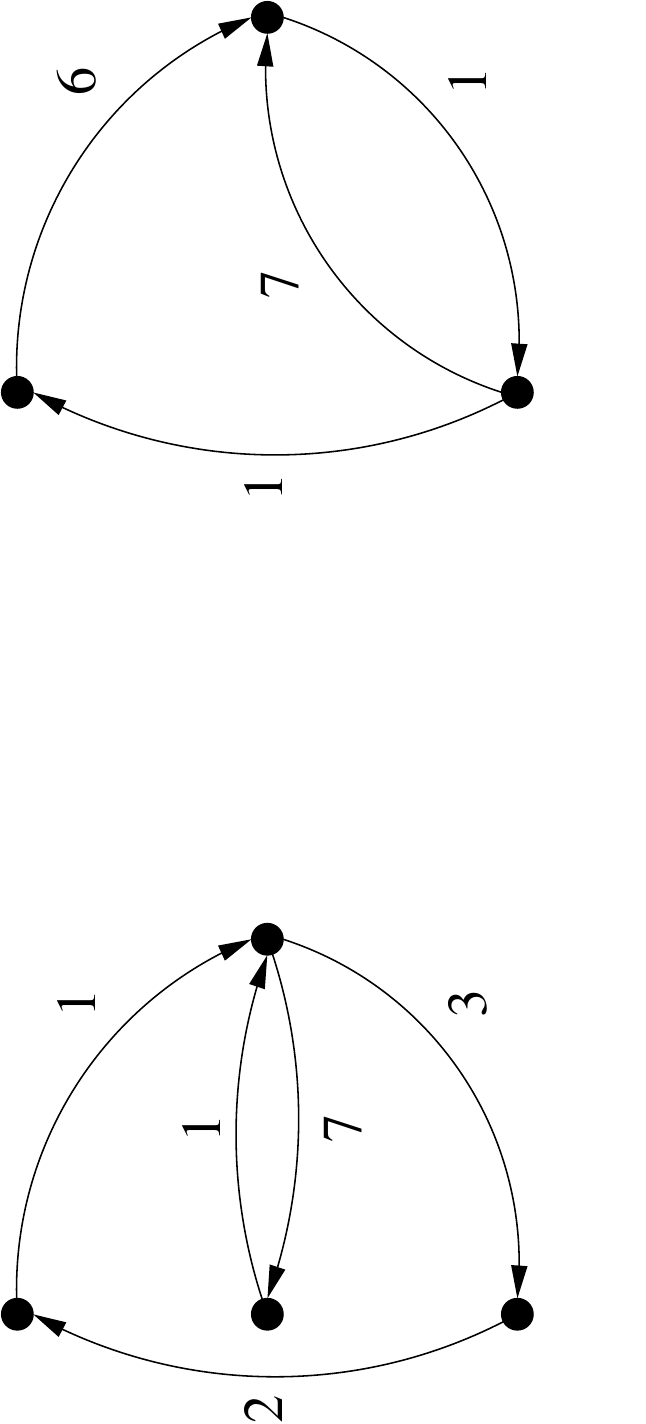}
    \put(19,3.5){$J$}
    \put(-3,-3){$\sigma(J)=\{3,-2,-1,0\}$}
    \put(82,3.5){$\mathcal{L}_S(J)$}
    \put(60,-3){$\sigma(\mathcal{L}_S(J))=\{3,-2,-1\}$}
    \put(37,25){$v_1$}
    \put(102,25){$v_1$}
    \end{overpic} 
\caption{For $S=\{v_1\}$ the graph $\mathcal{L}_S(J)$ is a reduction over the weight set of positive integers. Here only the structural set $S=\{v_1\}$ is labeled.}
\end{figure}

\textit{Step 3:} If the paths $P_1,\dots,P_k$ of $\pi\mathcal{B}_{ij}(G;S)$ are of the same length then in this construction there will be multiple edges between $v_i$ and these path's first interior vertices. In this case reduce these multiple edges to a single edge of weight $\sum_{i=1}^k\pi(P_i)$. The resulting graph is $\mathcal{L}_S(G)$. We note that the spectra $\sigma(G)$ and $\sigma(\mathcal{L}_S(G))$ are the same aside from a number of zeros. 

Now, suppose $U$ is a so called \textit{unital ring} \cite{7}, i.e. $U$ is a collection of numbers that contain all products of these numbers and sums of these products as well as the number (unit) 1. If $U$ is a unital ring containing all edge weights of $G$ then each edge weight of $\mathcal{L}_S(G)$ is in $U$. In this case we say $\mathcal{L}_S(G)$ is a reduction of $G$ \textit{over the weight set} $U$ if $\mathcal{L}_S(G)$ contains fewer vertices than $G$. As an example consider the graphs $J$  and $\mathcal{L}_S(J)$ in Fig. 4. 

For applications, the most important reductions over fixed weight sets are the cases in which $U$ is the ring of all positive numbers or the trivial ring consisting of the weight 1. When only the topology of a network is known but not the weights, it is natural to assign each edge the weight 1. In this case, a reduced network is again unweighted as each edge will have weight 1 (see Fig. 3 where $K$ is the reduction of $\mathcal{X}_S(K)$ over $U=\{1\}$). 

Isospectral expansion over fixed weight sets always exist if a network is not too simple (see Fact 5). This allows one to make the network sparser while maintaining the network's spectrum up to a known set. 

\begin{fact} \textbf{(Sparsification of Dynamical Networks)}
If a network has two cycles (closed paths) which are not loops but intersect in at least one vertex then such networks can be made (sequentially) more sparse. This is done via an isospectral expansion that exactly maintains the network's set of edge weights and preserves the network's spectrum up to a set of eigenvalues known in advance (see theorem 4.2 in \cite{3}). 
\end{fact}

This process is described as follows. For two paths $\beta$ and $\gamma$ in $\mathcal{B}_S(G)$ we say these paths are \textit{independent} if they share no interior vertices. An isospectral expansion of a graph $G$ over its weight set is the graph $\mathcal{X}_S(G)$ in which any two paths of $\mathcal{B}_S(G)$ have been made independent. That is, to each path in $\mathcal{B}_{ij}(G;S)$ there is a corresponding path in $\mathcal{B}_{ij}(\mathcal{X}_S(G);S)$ of the same length and edge weights. 

To describe the extent to which the eigenvalues of a graph differ from its expansion let $n_i$ be the number of paths in $\mathcal{B}_S(G)$ containing $v_i$. Then $\sigma(G)$ and $\sigma(\mathcal{X}_S(G))$ differ by $n_i-1$ eigenvalues of the form $\omega(e_{ii})$ for all $v_i\in\bar{S}$. As an example, in Fig. 3 we show an isospectral expansion $\mathcal{X}_S(K)$ of the graph $K$ over its weight set $\{1\}$. Here, each additional vertex in the expansion $\mathcal{X}_S(G)$ corresponds to an extra $0$ in its spectrum. Because $\mathcal{X}_S(G)$ is far more sparse than $K$ we refer to $\mathcal{X}_S(G)$ as a \textit{sparsification} of $K$. It is well known that there are various methods to analyze sparse networks. Our procedure allows us to make any network sparse.

In conclusion, we have presented a rigorous procedure which allows one to transform a network to a simpler one with fewer nodes and edges while preserving the network's collection of eigenvalues (spectrum). Such reductions can be used to establish spectral equivalence between networks and can be carried out according to any criterion related to any network characteristic (degree, betweenness, etc.). This allows one to compare the reductions of different networks and find similarities in their structure related to the chosen characteristics. This analysis is moreover facilitated by the fact that isospectral reductions visually simplify networks. These results can be readily applied to any real network. Additionally, this procedure can be used on dynamical networks to obtain improved stability results (see theorem 6.8 in \cite{3}).
 
The authors would like to thank M. Porter, C. Kirst, and M. Timme for useful suggestions. This work was partially supported by the NSF grant DMS-0900945 and the Humboldt Foundation.


\begin{thebibliography}{9}

\bibitem{1} R. Albert and A.-L. Barab\'{a}si, Rev. Mod. Phys. \textbf{74}, 47, (2002); S. Boccaletti, V. Latora, Y. Moreno, M. Chavez, and D. Hwang, Physics Reports \textbf{424}, 175.

\bibitem{5} S. N. Dorogovtsev and J. F. F. Mendes, Evolution of Networks: From Biological Nets to the Internet and WWW (Oxford University Press, USA, 2003);

\bibitem{4} M. E. J. Newman, A. L. Barab\'{a}si, and D. J. Watts, eds., The Structure and Dynamics of Networks  (Princeton University Press, 2006); M. E. J. Newman, Networks: An Introduction (Oxford University Press, 2010).

\bibitem{2} V. Afriamovich and L. Bunimovich, Nonlinearity (2007) \textbf{20} 1761; M. Blank and L. Bunimovich, Nonlinearity \textbf{19}, 329; J. G. Restrepo, E. Ott, and B. R. Hunt, Physical Review E \textbf{71}, 036151 , (2005).

\bibitem{8} T. Nishikawa1 and A. E. Motter, Phys. Rev. E \textbf{73}, 065106 (2006); T. Nishikawa1 and A. E. Motter, Physica D \textbf{224}, (2006).

\bibitem{3} L. Bunimovich and B. Webb, arXiv:1010.3272v1;

\bibitem{7} P. Cohn, An Introduction to Ring Theory (Springer, 2000).

\end{thebibliography}
\end{document}